\documentclass{amsproc}
\usepackage{graphicx}
\usepackage{amscd}
\usepackage{amsmath}
\usepackage{amsfonts}
\usepackage{amssymb}
\newtheorem{theorem}{Theorem}[section]
\newtheorem{corollary}[theorem]{Corollary}
\newtheorem{lemma}[theorem]{Lemma}
\newtheorem{proposition}[theorem]{Proposition}
\theoremstyle{definition}
\newtheorem{definition}[theorem]{Definition}

\theoremstyle{remark}
\newtheorem{remark}[theorem]{Remark}

\numberwithin{equation}{section}

\begin{document}
\title[Lefschetz Coincidence Theory]{Lefschetz Coincidence Theory for Maps Between Spaces of Different Dimensions}
\author{Peter Saveliev}
\address{Department of Mathematics, University of Illinois, 1409 West Green Street,
Urbana, IL 61801}
\email{saveliev@member.ams.org}

\begin{abstract}
For a given pair of maps $f,g:X\rightarrow M$ from an arbitrary topological
space to an $n$-manifold, the Lefschetz homomorphism is a certain graded
homomorphism $\Lambda_{fg}:H(X)\rightarrow H(M)$ of degree $(-n)$. We prove a
Lefschetz-type coincidence theorem: if the Lefschetz homomorphism is
nontrivial then there is an $x\in X$ such that $f(x)=g(x)$.
\end{abstract}\maketitle

\section{Introduction.}

Consider the\textbf{\ }Fixed Point Problem: ``If $X$ is a topological space
and $g:X\rightarrow X$ is a map, what can be said about the set $Fix(g)$ of
points $x\in X$ such that $g(x)=x$?'' The Coincidence Problem is concerned
with the same question about two maps $f,g:X\rightarrow Y$ and the set
$Coin(f,g)$ of $x\in X$ such that $f(x)=g(x).$

If $X$ is a sufficiently ``nice'' space (e.g., a polyhedron) then one may
associate to $g:X\rightarrow X$ an integer $\Lambda_{g},$ called the Lefschetz
number (see \cite{Brown}):
\[
\Lambda_{g}=L(g_{\ast})=\sum_{n}(-1)^{n}Trace(g_{\ast n}),
\]
where $g_{\ast n}$ is the endomorphism of the $n$-th homology group $H_{n}(X)
$ of $X$ induced by $f.$ Then the famous Lefschetz fixed point theorem states
that $\lambda_{g}\neq0\Rightarrow Fix(g)\neq\emptyset$. Now suppose we are
given a pair of continuous maps $f,g:X\longrightarrow Y,$ where only $Y$ has
to be a ``nice'' space and $X$ is arbitrary. Then one can define the Lefschetz
number $\Lambda_{fg}$ of the pair $(f,g)$ (see \cite[Section VI.14]{Bredon}):
\[
\Lambda_{fg}=L(g_{\ast}f_{!})=\sum_{n}(-1)^{n}Trace(g_{\ast n}f_{!}),
\]
where $f_{!}:H(Y)\rightarrow H(X)$ is a certain ``transfer'' homomorphism of
$f.$ Then a Lefschetz-type coincidence theorem states that $\lambda_{fg}%
\neq0\Rightarrow Coin(f,g)\neq\emptyset$.

Lefschetz coincidence theory has been developed for the following settings.

\begin{description}
\item [\textbf{Case 1}]$f:(M_{1},\partial M_{1})\rightarrow(M_{2},\partial
M_{2})$ is a boundary-preserving map between two $n$-manifolds with (possibly
empty) boundaries $\partial M_{1}$\ and $\partial M_{2}$, $g:M_{1}\rightarrow
M_{2}$ is arbitrary.
\end{description}

For closed manifolds, this is the setting of the original Lefschetz's result.
After many years, his theorem was extended to the case of manifolds with
boundary by M. Nakaoka \cite{Naka} and V. P. Davidyan \cite{David0},
\cite{David}.

\begin{description}
\item [\textbf{Case 2}]$f:X\rightarrow V$ maps a topological space to an open
subset of $\mathbf{R}^{n}$\ and all fibres $f^{-1}(y)$\ are acyclic,
$g:X\rightarrow V$ is compact.
\end{description}

This approach was developed by S. Eilenberg and D. Montgomery \cite{EM} and
later by L. Gorniewicz \cite{Gorn}, A. Granas \cite{GG2} and others. These
results treat fixed points of a multivalued map $G:Y\rightarrow Y$ by letting
$X$ be the graph of $G$ and $f,g$ be the projections, then $Fix(G)\equiv
\{x:x\in G(x)\}=Coin(f,g).$

In \cite{Sav} we proved a Lefschetz-type coincidence theorem that contains
Cases 1 and 2 and gave examples of coincidence situations not covered by the
known results (see \cite[Section 5]{Sav}). In particular we showed that the
projection of the torus $\mathbf{T}^{2}$ on the circle $\mathbf{S}^{1}$ has a
coincidence with any homologically trivial (inessential) map. This is an
example of a map between manifolds of different dimensions. In the present
paper we generalize the main results of \cite{Sav} in order to include a
Lefschetz-type coincidence theorem for the following setting.

\begin{description}
\item [\textbf{Case 3}]$f:T\times Y\rightarrow Y$\ is the projection, $Y$\ is
an ANR, $T$\ is any normal space, $g:T\times Y\rightarrow Y$ is arbitrary.
\end{description}

Here the coincidence set of the pair $(f,g),$ $Coin(f,g),$ is the fixed point
set $Fix(g)\equiv\{(t,x):g(t,x)=x\}$ of the ``parametrized'' map $g$. This
situation was studied by Knill \cite{Knill} and later by Geoghegan and Nicas
\cite{GN1,GN}, Geoghegan, Nicas and Oprea \cite{GNO}. The Lefschetz number is
replaced with a certain homomorphism $L(g):H(T)\rightarrow H(Y)$ of degree $0$
which is proven to be equal to the following.

\begin{definition}
\cite{Knill,GNO} The \textit{Knill trace} of $g$ is defined by
\[
L(g)(u)=\sum_{k\geq0}(-1)^{k+m}\sum_{j=1}^{\beta_{k}}x_{j}^{k}\frown g_{\ast
}(u\times b_{j}^{k}),
\]
where $u\in H_{m}(T)$ and for each $k\geq0$ $\{b_{j}^{k}:j=1,...,\beta_{k}\}$
is a basis for $H_{k}(Y)$ with corresponding dual basis $\{x_{j}%
^{k}:j=1,...,\beta_{k}\}$ for $H^{k}(Y).$ [Comment: The above definition uses
Spanier's sign convention \cite{Spanier}, we use Dold's sign convention
\cite{Dold} instead.]
\end{definition}

In this paper we use the results and techniques of our previous paper
\cite{Sav} to extend some of the definitions and theorems of \cite{GNO} to the
general case of two arbitrary maps $f,g:X\rightarrow Y$, i.e., $f$ is not
necessarily the projection, with the following reservation. For the sake of
simplicity, we limit our attention to the case when $Y=M$ is a manifold. Then
the Lefschetz homomorphism is a certain graded homomorphism $\Lambda
_{fg}:H(X,A)\rightarrow H(M)$ of degree $(-n)$, where $n=\dim M$. All relevant
examples assume that $X$ is a manifold as well and that
\[
\dim X>\dim M.
\]
We should tell from the start that the Lefschetz homomorphism of any pair of
maps $f,g:\mathbf{S}^{N}\rightarrow\mathbf{S}^{n}$ is trivial if $N\neq n.$

The author would like to thank Ross Geoghegan for a suggestion that lead to
this investigation. The author also thank the referee for a number of
suggestions that lead to significant improvements of the paper.

\textit{The Setup.} Throughout the paper we assume the following. By $H$
($H^{\ast}$) we denote the singular (co)homology with coefficients in a field
$R,$\ $M$ is an oriented connected compact $n$-manifold, $n\geq0$, with
boundary $\partial M$ and interior $\overset{\circ}{M},$ $O_{M}\in
H_{n}(M,\partial M)$ is the fundamental class of $(M,\partial M),$ $X$ is a
topological space, $A\subset X.$ We will study the existence of coincidences
of maps
\[
f:(X,A)\rightarrow(M,\partial M),\ g:X\rightarrow M,
\]
that satisfy
\[
Coin(f,g)\cap A=\emptyset.
\]

\section{The Lefschetz Class of a Homomorphism of Degree $m$.\label{LefClass}}

The following is a setup of the theory of the Lefschetz class of a graded
homomorphism of arbitrary degree as presented in \cite[Sections 1 and 3]{GNO}
(see also \cite{Sav} and \cite[p. 207-208]{Dold}). Let $E$ and $C$ be graded
$R$-spaces, $E$ finitely generated. Let $E^{\ast}$ denote the dual graded $R
$-space:
\[%
\begin{array}
[c]{l}%
E^{q}=Hom(E_{q}),\text{\quad}E^{\ast}=\{E^{q}\},\text{ and let}\\
(E^{\ast}\otimes E)_{m}=%
{\displaystyle\bigotimes}
{}_{q-p=m}(E^{p}\otimes E_{q}),\text{\quad}E^{\ast}\otimes E=\{(E^{\ast
}\otimes E)_{m}\}.
\end{array}
\]
Suppose we have two graded homomorphisms
\[%
\begin{array}
[c]{l}%
\frown:(E^{\ast}\otimes E)_{m}\longrightarrow C_{m}\text{, and}\\
\theta:(E^{\ast}\otimes E)_{m}\longrightarrow Hom_{m}(E,E)
\end{array}
\]
where $Hom_{m}(E,E)$ denotes the space of all graded homomorphisms of degree
$m$. We define $\theta$ as follows:
\[
\lbrack\theta(a\otimes b)](u)=(-1)^{\mid b\mid\cdot\mid u\mid}a(u)\cdot b,
\]
where $a\in E^{k},b\in E_{m+k},u\in E_{k},$ $a\otimes b\in(E^{\ast}\otimes
E)_{m},$ $|w|$ stands for the degree of $w$. Then by \cite[Proposition
VII.6.3, p. 208]{Dold}, $\theta$ is an isomorphism.

\begin{definition}
\label{DefLefClass}For an endomorphism $h:E\longrightarrow E$ of degree $m$ of
a finitely generated graded module $E$, we define the\textit{\ Lefschetz
class} of $h$ by
\[
L(h)=\frown\theta^{-1}(h)\in C_{m}.
\]
\end{definition}

The following representation of the Lefschetz class as a Knill-like trace is
proven similarly to \cite[Proposition 1.2]{GNO}.

\begin{proposition}
\label{KnillTrace}Let $h:E\rightarrow E$ be a homomorphism of degree $m.$ Let
$\{a_{1}^{k},...,a_{m_{k}}^{k}\}$ be a basis for $E_{k}$ and $\{x_{1}%
^{k},...,x_{m_{k}}^{k}\}$ the corresponding dual basis for $E^{k}.$ Then
\[
L(h)=\sum_{k}(-1)^{k(k+m)}\sum_{j}x_{j}^{k}\frown h(a_{j}^{k}).
\]
\end{proposition}

When $m=0,$ $C_{m}=R$ and $\frown=e$ is the evaluation map, we have the usual
Lefschetz number as the alternating sum of traces, see \cite[p. 208]{Dold}.

\section{The Evaluation Formula for the Lefschetz Class.\label{EvalForm}}

For any space $Y$ we define the following functions:

\noindent\textit{the transposition} $t:Y\times Y\longrightarrow Y\times Y$
given by $t(x,y)=(y,x);$

\noindent\textit{the diagonal map} $\delta:Y\rightarrow Y\times Y$ given by
$\delta(x)=(x,x);$

\noindent\textit{the scalar multiplication} $q:R\mathbf{\otimes}%
H(Y)\longrightarrow H(Y)$ given by $q(r\otimes v)=r\cdot v;$

\noindent\textit{the tensor multiplication }$O_{M}^{\times}%
:H(Y)\longrightarrow H(M,\partial M)\otimes H(Y)$ given by $O_{M}^{\times
}(v)=O_{M}\otimes v;$

\noindent\textit{the projection} $P_{k}:H(Y)\longrightarrow H_{k}(Y)$, $k\geq0.$

Let
\[
M^{\prime}=M\cup C,
\]
where $C=\partial M\times\lbrack0,1)$ is the collar attached to the boundary
of $M.$ Define
\[
M^{\times}=(M\times M^{\prime},M\times M^{\prime}\backslash\delta(M^{\prime})
\]
and \noindent\textit{the inclusions} $i:\overset{\circ}{M}\longrightarrow
M,\ I:(M,\partial M)\times\overset{\circ}{M}\longrightarrow M^{\times}.$ If
$\pi:M\times M^{\prime}\rightarrow M$ is the projection on the first factor
then $\zeta=(M,\pi,M\times M^{\prime},\delta)$ is the tangent microbundle of
$M$ \cite[Chapter 14]{Switzer} and \noindent\textit{the Thom isomorphism}
$\varphi:H(M^{\times})\rightarrow H(M)$ is given by $\varphi(x)=\pi_{\ast
}(\tau\frown x),$ where and $\tau$ is the Thom class of $\zeta.$

The proofs in this section follow the ideas of Dold \cite{Dold0} and
\cite[Section VI.6]{Dold}, see also \cite{Gorn} and \cite{Sav}.

\begin{lemma}
[Generalized Dold's Lemma]\label{MDold}Suppose that the map $\Phi
:H(M)\longrightarrow H(M)$ is given as the composition of the following
homomorphisms:
\[%
\begin{array}
[c]{l}%
\Phi:H(\overset{\circ}{M})^{\underrightarrow{~O_{M}^{\times}~~}}H(M,\partial
M)\otimes H(\overset{\circ}{M})^{\underrightarrow{~\delta_{\ast}\otimes Id~~}%
}H(M,\partial M)\otimes H(M)\otimes H(\overset{\circ}{M})\\
^{\underrightarrow{~Id\otimes t_{\ast}~~}}H(M,\partial M)\otimes
H(\overset{\circ}{M})\otimes H(M)^{\underrightarrow{~I_{\ast}\otimes Id~~}%
}H(M^{\times})\otimes H(M)\\
^{\underrightarrow{~P_{n}\otimes Id~~}}H_{n}(M^{\times})\otimes H(M)\simeq
R\otimes H(M)\overset{q}{\longrightarrow}H(M).
\end{array}
\]
Then
\[
\Phi=i_{\ast}.
\]
\end{lemma}

\begin{proof}
In \cite[Theorem 7.2]{Sav} we proved this statement for the following
composition:
\[%
\begin{array}
[c]{l}%
\Psi:H(K)^{\underrightarrow{~O_{K}^{\times}~~}}H(V,V\backslash K)\otimes
H(K)^{\underrightarrow{~\delta_{\ast}\otimes Id~~}}H(V,V\backslash K)\otimes
H(V)\otimes H(K)\\
^{\underrightarrow{~Id\otimes t_{\ast}~~}}H(V,V\backslash K)\otimes
H(K)\otimes H(V)^{\underrightarrow{~I_{\ast}\otimes Id~~}}H(M^{\times})\otimes
H(V)\\
^{\underrightarrow{~P_{n}\otimes Id~~}}H_{n}(M^{\times})\otimes H(V)\simeq
R\otimes H(V)\overset{q}{\longrightarrow}H(V),
\end{array}
\]
where $M$ is a closed manifold and $(M,V,V\backslash K)$ is an excisive triad.
If $M$ has the boundary the proof is the same except for the last step: for
any $p\in\overset{\circ}{M},$ $I_{\ast}(O_{M}\otimes p)=1.$ This formula
follows from definitions of the Thom class and the fundamental class, see
\cite[Chapter 14]{Switzer}. Now we obtain the above statement by putting
$V=M,K=\overset{\circ}{M}.$
\end{proof}

Now let $E=C=H(M)$ and let $\frown$ be the usual cap-product.

For the following three lemmas we fix $m\geq0$ and let $h:H_{\ast
}(M)\rightarrow H_{\ast+m}(M)$ be a homomorphism of degree $m$.

\begin{lemma}
\label{diag3}There is a homomorphism $J:H_{n-\ast}(M,\partial
M)\longrightarrow H^{\ast}(M)$ such that the following diagram commutes
\[%
\begin{array}
[c]{ccc}%
H_{n-\ast}(M,\partial M)\otimes H_{m+\ast}(\overset{\circ}{M}) &
^{\underrightarrow{\qquad J\otimes i_{\ast}\qquad}} & H^{\ast}(M)\otimes
H_{m+\ast}(M)\\
\downarrow^{I_{\ast}} &  & \downarrow^{\frown}\\
H_{n+m}(M^{\times}) & ^{\underrightarrow{\qquad\quad\varphi\quad\qquad}} &
H_{m}(M),
\end{array}
\]
or
\[
\frown(J\otimes i_{\ast})=\varphi I_{\ast}.
\]
\end{lemma}

\begin{proof}
Let $J$ be given by
\[
J(u)(v)=\varphi I_{\ast}(u\otimes v),\text{\quad}u\in H_{n-j}(M,\partial
M),v\in H_{m+j}(M).
\]
Let's consider the diagram on the chain level and start in the left upper
corner with $u\otimes v,$ where $u$ is a $(n-j)$-chain and $v$ is a
$(m+j)$-chain. Then, in terms of the Alexander-Whitney approximation, we get
$\pi_{\ast}\tau I_{\ast}(u\otimes_{j}\lfloor v)v\rfloor_{m}$ in the right
lower corner.
\end{proof}

\begin{lemma}
\label{defa}Let
\[
a=(J\otimes h)\delta_{\ast}(O_{M})\in(H^{\ast}(M)\otimes H(M))_{m}.
\]
Then
\[
\frown(a)=\varphi I_{\ast}(Id\otimes i_{\ast}^{-1}h)\delta_{\ast}(O_{M})\in
H_{m}(M).
\]
\end{lemma}

\begin{proof}
By Lemma \ref{diag3} we have
\[
\frown(a)=\frown(J\otimes i_{\ast})(Id\otimes i_{\ast}^{-1}h)\delta_{\ast
}(O_{M})=\varphi I_{\ast}(Id\otimes i_{\ast}^{-1}h)\delta_{\ast}(O_{M}).
\]
\end{proof}

\begin{lemma}
\label{lastdiag}The following diagram commutes:
\[%
\begin{array}
[c]{ccc}%
H_{n-k}(M,\partial M)\otimes H_{k}(M)\otimes H_{p}(\overset{\circ}{M}) &
^{\underrightarrow{J\otimes h\otimes i_{\ast}}} & H^{k}(M)\otimes
H_{m+k}(M)\otimes H_{p}(M)\\
~~~~~~~\downarrow^{Id\otimes t_{\ast}} &  & ~~~~~~~~\downarrow^{Id\otimes
t_{\ast}}\\
H_{n-k}(M,\partial M)\otimes H_{p}(\overset{\circ}{M})\otimes H_{k}(M) &
^{\underrightarrow{J\otimes i_{\ast}\otimes h}} & H^{k}(M)\otimes
H_{p}(M)\otimes H_{m+k}(M)\\
~~~~~~~~~~~~~~~~~~~~\downarrow^{I_{\ast}\otimes Id} &  & ~~~~~~~~~~~\downarrow
^{\frown\otimes Id}\\
H_{n-k+p}(M^{\times})\otimes H_{k}(M) & ^{\underrightarrow{~\varphi\otimes
h}~} & H_{p-k}(M)\otimes H_{m+k}(M)\\
~~~~~\downarrow^{q(P_{n}\otimes Id)} &  & ~~~~~\downarrow^{q(P_{0}\otimes
Id)}\\
H_{k}(M) & ^{\underrightarrow{\quad h\quad}} & H_{m+k}(M).
\end{array}
\]
\end{lemma}

\begin{proof}
The first square trivially commutes. For the second, going $\downarrow
\rightarrow,$ we get $\varphi I_{\ast}\otimes h.$ Now going $\rightarrow
\downarrow$, from Lemma $\ref{diag3}$ we get$:$%
\[
(\frown\otimes Id)(J\otimes i_{\ast}\otimes h)=(\frown(J\otimes i_{\ast
})\otimes h)=\varphi I_{\ast}\otimes h.
\]
The third square commutes as $p=k.$
\end{proof}

\begin{theorem}
[Evaluation Formula]\label{EvalFormula}For any homomorphism $h:H(M)\rightarrow
H(M)$ (of degree $m$) we have
\[
L(h)=\varphi I_{\ast}(Id\otimes i_{\ast}^{-1}h)\delta_{\ast}(O_{M}).
\]
\end{theorem}

\begin{proof}
(cf. \cite[Theorem 8.5]{Sav}) We start in the left upper corner of the above
diagram with $\delta_{\ast}(O_{M})\otimes v$, where $v\in H(M).$ Let
$u=i_{\ast}(v).$ Then going $\downarrow\rightarrow,$ we get $h(u)$ by Lemma
\ref{MDold}. Let $a=%
{\textstyle\sum\nolimits_{i}}
a_{i-m}\otimes a_{i}^{\prime}$ with $a_{i-m}\in H^{i-m}(M),a_{i}^{\prime}\in
H_{i}(M).$ Then going $\rightarrow\downarrow,$ we get:
\[%
\begin{array}
[c]{l}%
\quad\ q(P_{0}\otimes Id)(\frown\otimes Id)(Id\otimes t_{\ast})(J\otimes
h\otimes Id)(\delta_{\ast}(O_{M})\otimes u)\\
=q(P_{0}\frown\otimes Id)(Id\otimes t_{\ast})((J\otimes h)\delta_{\ast}%
(O_{M})\otimes u)\\
=q(P_{0}\frown\otimes Id)(Id\otimes t_{\ast})(a\otimes u)\qquad\text{by
definition of }a\text{ (Lemma }\ref{defa})\\
=q(P_{0}\frown\otimes Id)(Id\otimes t_{\ast})(%
{\textstyle\sum\nolimits_{i}}
a_{i-m}\otimes a_{i}^{\prime}\otimes u)\\
=q(P_{0}\frown\otimes Id)%
{\textstyle\sum\nolimits_{i}}
(-1)^{\mid a_{i}^{\prime}\mid\cdot\mid u\mid}(a_{i-m}\otimes u\otimes
a_{i}^{\prime})\qquad\\
=%
{\textstyle\sum\nolimits_{i}}
(-1)^{\mid a_{i}^{\prime}\mid\cdot\mid u\mid}P_{0}(a_{i-m}\frown u)\cdot
a_{i}^{\prime}\qquad\qquad\\
=%
{\textstyle\sum\nolimits_{i}}
(-1)^{\mid a_{i}^{\prime}\mid\cdot\mid u\mid}a_{i-m}(u)\cdot a_{i}^{\prime}\\
=%
{\textstyle\sum\nolimits_{i}}
\theta(a_{i-m}\otimes a_{i}^{\prime})(u)\qquad\qquad\\
=\theta(a)(u).
\end{array}
\]
Thus $\theta(a)=h:H(M)\longrightarrow H(M).$ Therefore by definition of the
Lefschetz homomorphism we have
\[
L(h)=L(\theta(a))=\frown(a).
\]
Now the statement follows from Lemma \ref{defa}.
\end{proof}

\section{Transfers and the Coincidence Homomorphism.\label{Transfers}}

Since $Coin(f,g)\cap A=\emptyset,$ the map $(f\times g)\delta
:(X,A)\longrightarrow M^{\times}$ is well defined.

\begin{definition}
\label{DefCoin}The \textit{coincidence homomorphism} $I_{fg}$ \textit{of the
pair} $(f,g)$ is the homomorphism $I_{fg}:H(X,A)\rightarrow H(M^{\times})$ of
degree $0$ defined by
\[
I_{fg}=(f\times g)_{\ast}\delta_{\ast}.
\]
\end{definition}

It is clear that $I_{fg}\neq0\Longrightarrow Coin(f,g)\neq\emptyset.$

From dimensional considerations we obtain the following.

\begin{proposition}
\label{Trivial1}If $z\in H_{i}(X,A)$, $i<n$ or $i>2n,$ then $I_{fg}(z)=0.$
\end{proposition}

This is the reason why it is often sufficient to consider \textit{the
coincidence index with respect to} $\mu\in H_{n}(X,A)$, as in \cite{Sav}:
\[
I_{fg}(\mu)=(f\times g)_{\ast}\delta_{\ast}(\mu).
\]
In fact when $X$ is a manifold and $\deg X=n=\deg M$, the best we can do is to
consider $I_{fg}(O_{X})$, where $O_{X}$ is the fundamental class of $X$ (see
also Corollaries \ref{Trivial3} and \ref{Trivial4})$.$

\begin{corollary}
$I_{fg}$ is trivial for maps $\mathbf{S}^{N}\rightarrow M$ if $N>2n.$
\end{corollary}

\begin{definition}
\label{DefTrans}The \textit{transfer (or the shriek map) of }$f$\textit{\ with
respect to} $z\in H_{n+m}(X,A)$ is the homomorphism $f_{!}^{z}:H_{\ast
}(M)\longrightarrow H_{\ast+m}(X)$ of degree $m$ given by
\[
f_{!}^{z}=(f^{\ast}D^{-1})\frown z,
\]
where $D:H^{\ast}(M,\partial M)\rightarrow H_{n-\ast}(M)$ is the
Poincare-Lefschetz duality isomorphism $D(x)=x\frown O_{M}$.
\end{definition}

When $X$ is a manifold and $z$ is its fundamental class, $f_{!}^{z}=f_{!}$ is
the usual transfer homomorphism, or an \textit{Umkehr}-homomorphism, of $f $
\cite[Section VIII.10]{Dold}, or a shriek map \cite[p. 368]{Bredon}.

We can represent the coincidence homomorphism via the transfer as follows.

\begin{theorem}
[Representation Formula]\label{representation}
\[
I_{fg}(z)=I_{\ast}(Id\otimes i_{\ast}^{-1}g_{\ast}f_{!}^{z})\delta_{\ast
}(O_{M}).
\]
\end{theorem}

\begin{proof}
In \cite[Theorem 2.1]{Sav} the formula is proven for $z\in H_{n}(X,A),$ but
the proof is valid for any $z\in H(X,A).$
\end{proof}

Thus $I_{fg}(z)$ is the image of $O_{M}$ under the composition of the
following maps:
\[
H(M,\partial M)^{\underrightarrow{\ \delta_{\ast}\ }}H(M,\partial M)\otimes
H(M)^{\underrightarrow{Id\otimes i_{\ast}^{-1}h\ \ }}H(M,\partial M)\otimes
H(\overset{\circ}{M})^{\underrightarrow{\ I_{\ast}\ }}H(M^{\times}),
\]
where homomorphism $h=g_{\ast}f_{!}^{z}:H(M)\rightarrow H(M)$ of degree $m$ is
defined by the following diagram:
\[%
\begin{array}
[c]{ccc}%
H^{\ast}(X,A) & ^{\ \underleftarrow{\quad\ \ \ \ \ f^{\ast}\ \quad}}\  &
H^{\ast}(M,\partial M)\\
~~~~\downarrow^{\frown z} &  & \downarrow^{D}\\
H(X) & ^{\underrightarrow{\ \ \ \ \ \ g_{\ast}\ \ \ \ \ \ }} & H(M).
\end{array}
\]

For a given map $f:X\rightarrow Y$, Gottlieb \cite{Gottlieb} (see also
\cite{BG} and \cite{Gottlieb1}) defines \textit{a partial transfer of }%
$f$\textit{\ with trace }$k$ as a homomorphism $\tau:H(Y)\rightarrow H(X)$
such that $f_{\ast}\tau:H(Y)\rightarrow H(Y)$ is the multiplication by $k$:
\[
f_{\ast}\tau=k\cdot Id.
\]
Then independently of Theorem \ref{representation} we can prove its analogue
(for $m=0$):

\begin{proposition}
Suppose $\tau$ is a partial transfer of $f:X\rightarrow M$ ($\partial
M=\emptyset$) with trace $k\neq0.$ Then
\[
I_{fg}(z)=I_{\ast}(Id\otimes g_{\ast}\tau)\delta_{\ast}(O_{M}),
\]
where $z=\tau(O_{M}).$
\end{proposition}

\begin{proof}
Consider the following commutative diagram:
\[%
\begin{tabular}
[c]{lllllll}%
&  & $H(M)$ & $^{\underrightarrow{\ \ \ \ \delta_{\ast}\ \ \ \ }}$ &
$H(M)\otimes H(M)$ &  & \\
& $\nearrow_{k\cdot Id}$ & $\uparrow_{f_{\ast}}$ &  & $\uparrow_{f_{\ast
}\otimes f_{\ast}}$ & $\nwarrow^{k^{2}\cdot Id}$ & \\
$H(M)$ & $^{\underrightarrow{\ \ \ \ \ \tau\ \ \ \ \ }}$ & $H(X)$ &
$^{\underrightarrow{\ \ \ \ \delta_{\ast}\ \ \ \ }}$ & $H(X)\otimes H(X)$ &
$^{\ \underleftarrow{\ \ \tau\otimes\tau\ }}$ & $H(M)\otimes H(M)$%
\end{tabular}
\]
Then going from the left to the right we get $\frac{1}{k}\delta_{\ast}.$ Hence
$\delta_{\ast}\tau=(\tau\otimes\tau)\frac{1}{k}\delta_{\ast},$ therefore
$k\cdot\delta_{\ast}\tau=(\tau\otimes\tau)\delta_{\ast}$. Thus the next
diagram is commutative:
\[%
\begin{array}
[c]{ccc}%
H(M) & ^{\underrightarrow{\ \ \ \ \ \tau\ \ \ \ \ }} & H(X)\\
\downarrow^{\delta_{\ast}} &  & \ \ \downarrow^{k\cdot\delta_{\ast}}\\
H(M)\otimes H(M) & ^{\underrightarrow{\ \ \ \tau\otimes\tau\ \ \ \ }} &
H(X)\otimes H(X)\\
& \searrow^{k\cdot Id\otimes\tau} & \qquad\downarrow^{f_{\ast}\otimes Id}\\
&  & H(M)\otimes H(X).
\end{array}
\]
Then
\[
(Id\otimes\tau)\delta_{\ast}=(f_{\ast}\otimes Id)\delta_{\ast}\tau,
\]
so
\[
I_{\ast}(Id\otimes g_{\ast}\tau)\delta_{\ast}=I_{\ast}(f_{\ast}\otimes
g_{\ast})\delta_{\ast}\tau.
\]
Now the statement follows from the definition of $I_{fg}$.
\end{proof}

The transfer $f_{!}^{z}$ of $f$ defined above can be an example of a partial
transfer (see \cite[Proposition 1]{Gottlieb} or \cite[Proposition VI.14.1 (6),
p. 394]{Bredon}):

\begin{proposition}
\label{f-inverse}If there is a $z\in H(X,A)$ such that $f_{\ast}(z)=k\cdot
O_{M}$ then $f_{!}^{z}$ is a partial transfer of $f$ with trace $k.$
\end{proposition}

Thus any partial transfer with nonzero trace satisfies the statement of
Theorem \ref{representation} and, therefore, the Lefschetz-type Theorem
\ref{LefThm} below. On the other hand, $k\neq0$ is a strong restriction as it
implies that $f_{*}(z)=k\cdot O_{M}$, which in case of manifolds of equal
dimensions means that $\deg f\neq0.$

\section{The Lefschetz Homomorphism of the Pair.\label{Manifolds}}

For a fixed $z\in H_{n+m}(X,A),$ the homomorphism $g_{\ast}f_{!}%
^{z}:H(M)\rightarrow H(M)$ has degree $m.$ Then by Definition
\ref{DefLefClass} we have $L(g_{\ast}f_{!}^{z})\in H_{m}(M).$

\begin{definition}
\label{LefHomo}The \textit{Lefschetz homomorphism }$\Lambda_{fg}:H_{\ast
}(X,A)\rightarrow H_{\ast-n}(M)$\textit{\ of the pair }$(f,g)$ is the
homomorphism of degree $(-n)$ given by
\[
\Lambda_{fg}(z)=L(g_{\ast}f_{!}^{z}),\ z\in H(X,A).
\]
\end{definition}

\begin{remark}
Suppose $N>n.$ Then from dimensional considerations it follows that
$\Lambda_{fg}=0$ for maps $f,g:\mathbf{S}^{N}\rightarrow\mathbf{S}^{n}.$
\end{remark}

\begin{remark}
\label{Trivial2}Suppose $z\in H_{i}(M,\partial M)$ and $i>n.$ Then $z=0$,
therefore $\Lambda_{Id,Id}(z)=0.$
\end{remark}

\begin{remark}
\label{KnillTrace2}Since the degree of $g_{\ast}f_{!}^{z}$ is $|z|-n,$ the
Lefschetz homomorphism is represented as a Knill-like trace (Proposition
\ref{KnillTrace}):
\[
\Lambda_{fg}(z)=\sum_{k}(-1)^{k(k+|z|-n)}\sum_{j}x_{j}^{k}\frown g_{\ast}%
f_{!}^{z}(a_{j}^{k}),
\]
where $\{a_{1}^{k},...,a_{m_{k}}^{k}\}$ is a basis for $H_{k}(M)$ and
$\{x_{1}^{k},...,x_{m_{k}}^{k}\}$ the corresponding dual basis for $H^{k}(M).$
\end{remark}

The Lefschetz homomorphism satisfies a naturality property below. The formula
is similar to the one in \cite[Theorem 2.6]{GNO} but the proof is much shorter
because we do not use the Knill trace (see also Theorem \ref{Natur2}).

\begin{theorem}
[Naturality I]\label{Natur}Let $(Y,B)$ be a topological space and
$h:(Y,B)\rightarrow(X,A) $ be a map. Then
\[
\Lambda_{fh,gh}=\Lambda_{fg}h_{\ast}.
\]
\end{theorem}

\begin{proof}
Let $z\in H(Y,B)$. Then we have
\[
(gh)_{\ast}(fh)_{!}^{z}=g_{\ast}h_{\ast}(h^{\ast}f^{\ast}D^{-1}\frown
z)=g_{\ast}(f^{\ast}D^{-1}\frown h_{\ast}(z))=g_{\ast}f_{!}^{h_{\ast}(z)}.
\]
\end{proof}

\begin{corollary}
\label{NaturTr}If $h$ has a partial transfer $\tau$ with trace $k$ then
\[
\Lambda_{fh,gh}\tau=k\cdot\Lambda_{fg}.
\]
\end{corollary}

Propositions \ref{Natur} and \ref{NaturTr} generalize the well known formula
for maps between two $n$-manifolds \cite[Corollary VI.14.6, p. 297]{Bredon}:
\[
L(fh,gh)=\deg(h)L(f,g),
\]
where $L(f,g)$ is the ordinary Lefschetz number.

The following corollary shows how the Lefschetz homomorphism generalizes the
Knill trace of a parametrized map defined in \cite[Definition 2.1]{GNO}, see
also Corollary \ref{Knill2}.

\begin{corollary}
\label{Knill1}Suppose $(X,A)=Y\times(M,\partial M),$ $g:X\rightarrow M$ is a
map and $p:Y\times(M,\partial M)\rightarrow(M,\partial M)$ is the projection
(then $Fix(g)=Coin(p,g)$). Then
\[
\Lambda_{pg}(u\times O_{M})=L(g_{u\ast}),\ u\in H(Y),
\]
where $g_{u}:H(M)\rightarrow H(M)$ is given by $g_{u}(x)=(-1)^{(n-|x|)|u|}%
g_{\ast}(u\times x).$
\end{corollary}

\begin{proof}
Let $x\in H(M)$ and suppose $x=z\frown O_{M}$ for some $z\in H^{\ast}(M).$
Then we have
\begin{align*}
p^{\ast}(z)  &  \frown(u\times O_{M})=(1\times z)\frown(u\times O_{M}%
)=(-1)^{|z||u|}(1\frown u)\times(z\frown O_{M})\\
&  =(-1)^{(n-|x|)|u|}u\times x.
\end{align*}
Therefore
\[
g_{\ast}p_{!}^{u\times O_{M}}(x)=g_{\ast}(p^{\ast}(z)\frown(u\times
O_{M}))=(-1)^{(n-|x|)|u|}g_{\ast}(u\times x)
\]
and the statement follows.
\end{proof}

\section{Further Properties.}

Theorems \ref{EvalFormula} and \ref{representation} imply the following theorem.

\begin{theorem}
[Lefschetz-Type Coincidence Theorem]\label{LefThm}The coincidence homomorphism
is equal to the Lefschetz homomorphism:
\[
\varphi I_{fg}=\Lambda_{fg}.
\]
Moreover, if $\Lambda_{fg}\neq0,$ then $(f,g)$ has a coincidence.
\end{theorem}

According to Corollary \ref{Knill1}, $\Lambda_{pg}$ generalizes the Knill
trace of a parametrized map $g:Y\times M\rightarrow M$, while in \cite{GNO}
the Knill trace is defined for $F:Y\times(X,A)\rightarrow(X,A)$, i.e., as a
map of pairs. But the result corresponding to the one above is due to Knill
\cite[Theorem 1]{Knill} and is proven for the case of $F:Y\times X\rightarrow
X$.

The above identity allows us to establish some facts about the Lefschetz
homomorphism that are hard to obtain directly from its definition.

\begin{theorem}
[Naturality II]\label{Natur2}(cf. \cite[Theorem 2.6]{GNO}) Let $(X^{\prime
},A^{\prime})$ be a topological space, $(M^{\prime},\partial M^{\prime})$
another $n$-manifold, $h:(X,A)\rightarrow(X^{\prime},A^{\prime}),$
$k:(M,\partial M)\rightarrow(M^{\prime},\partial M^{\prime}),$ $f^{\prime
},g^{\prime}:(X^{\prime},A^{\prime})\rightarrow(M^{\prime},\partial M^{\prime
})$ maps, and $f^{\prime}h=kf,$ $g^{\prime}h=kg$, i.e., we have the following
two (in one) commutative diagrams:
\[%
\begin{array}
[c]{ccc}%
(X,A) & ^{\underrightarrow{\ \ \ \ \ f,g\ \ \ \ \ }} & (M,\partial M)\\
\downarrow^{h} &  & \ \ \downarrow^{k}\\
(X^{\prime},A^{\prime}) & ^{\underrightarrow{\ \ \ f^{\prime},g^{\prime
}\ \ \ \ }} & (M^{\prime},\partial M^{\prime}).
\end{array}
\]
Suppose also that $k$ is a homeomorphism. Then
\[
k_{\ast}\Lambda_{fg}=\Lambda_{f^{\prime}g^{\prime}}h_{\ast}.
\]
\end{theorem}

\begin{proof}
The fact that $k$ is a homeomorphism implies two things. First, $k\times
k:M^{\times}\rightarrow(M^{\prime})^{\times}$ is well defined. Hence from the
naturality of the Thom isomorphism we have
\[
k_{\ast}\varphi=\varphi^{\prime}(k\times k)_{\ast},
\]
where $\varphi^{\prime}$ is the Thom isomorphism for $M^{\prime}.$ Second,
since $Coin(f,g)\cap A=\emptyset,$ it follows that $Coin(kf,kg)\cap
A=\emptyset$ and, therefore, $I_{kf,kg}$ is well defined (Definition
\ref{DefCoin}). In the computation below we also use Theorem \ref{Natur} (a
similar statement can be proven independently for $I_{fg}$), Theorem
\ref{LefThm} and the trivial fact that $I_{kf,kg}=(k\times k)_{\ast}I_{fg}.$
We have
\[
\Lambda_{f^{\prime}g^{\prime}}h_{\ast}=\Lambda_{f^{\prime}h,g^{\prime}%
h}=\Lambda_{kf,kg}=\varphi^{\prime}I_{kf,kg}=\varphi^{\prime}(k\times
k)_{\ast}I_{fg}=k_{\ast}\varphi I_{fg}=k_{\ast}\Lambda_{fg}.\text{ }%
\]
Observe that the Lefschetz homomorphism $\Lambda_{fg}$ is well defined without
the restriction $Coin(f,g)\cap A=\emptyset.$
\end{proof}

Even when $A=\partial M=\emptyset$, the definition of the Lefschetz
homomorphism is not symmetric, but the one of the coincidence homomorphism
(Definition \ref{DefCoin}) is, as follows:
\[
I_{fg}(z)=t_{\ast}I_{gf}(z).
\]
Now we use the fact that $t_{\ast}(x)=(-1)^{n}x$ for $x\in H(M^{\times})$ (the
proof of this formula is dual to the proof of Lemma 6.16 of \cite[p.
165]{Vick}). As a result we have the following property.

\begin{proposition}
[Symmetry]\label{Symmetry}Suppose $f,g:X\rightarrow M$ are maps ($\partial
M=\emptyset).$ Then
\[
\Lambda_{fg}(z)=(-1)^{n}\Lambda_{gf}(z),\ z\in H(X).
\]
\end{proposition}

It follows that $\Lambda_{ff}=0$ when $n$ is odd (in particular, $\chi(M)=0)$.

Now we can obtain another representation of the Knill trace of a parametrized
map, in terms of the Lefschetz homomorphism (cf. Corollary \ref{Knill1}):

\begin{corollary}
\label{Knill2}Suppose $g:Y\times M\rightarrow M$ ($\partial M=\emptyset)$ is a
map and $p:Y\times M\rightarrow M$ is the projection. Then
\[
L(g_{u\ast})=(-1)^{n}\Lambda_{gp}(u\times O_{M}),\ u\in H(Y),
\]
where $g_{u}:H(M)\rightarrow H(M)$ is given by $g_{u}(x)=(-1)^{(n-|x|)|u|}%
g_{\ast}(u\times x).$
\end{corollary}

\begin{proof}
It follows from Corollary \ref{Knill1} and Proposition \ref{Symmetry}.
\end{proof}

The proof of following property of the coincidence homomorphism is trivial.

\begin{theorem}
[Product Theorem](cf. \cite[Theorem 4.5]{GNO}) Let $(X^{\prime},A^{\prime})$
be a topological space, $(M^{\prime},\partial M^{\prime})$ a manifold,
$f^{\prime},g^{\prime}:(X^{\prime},A^{\prime})\rightarrow(M^{\prime},\partial
M^{\prime})$ maps. Then there is a commutative diagram:
\[%
\begin{array}
[c]{ccc}%
H((X,A)\times(X^{\prime},A^{\prime})) & ^{\underrightarrow
{\ \ \ \ \ I_{f\times f^{\prime},g\times g^{\prime}}\ \ \ \ \ }} & H((M\times
M^{\prime})^{\times})\\
\uparrow^{\xi} &  & \ \ \uparrow^{C_{\ast}\eta}\\
H(X,A)\otimes H(X^{\prime},A^{\prime}) & ^{\underrightarrow{\ \ \ \ \ I_{fg}%
\otimes I_{f^{\prime}g^{\prime}}\ \ \ \ \ \ }} & H(M^{\times})\otimes
H((M^{\prime})^{\times}),
\end{array}
\]
where $C:M^{\times}\times(M^{\prime})^{\times}\rightarrow(M\times M^{\prime
})^{\times}$ is the map which interchanges the middle factors, $\xi$ and
$\eta$ are the K\"{u}nneth isomorphisms.
\end{theorem}

Under certain circumstances the Lefschetz homomorphism is trivial in all
dimensions but $n.$

\begin{corollary}
\label{Trivial3}If $z\in H_{i}(X,A)$, $i\neq n,$ then $\Lambda_{ff}(z)=0.$
\end{corollary}

\begin{proof}
If $i>n$ then $\Lambda_{ff}(z)=\Lambda_{Id,Id}f_{\ast}(z)=0$ by Proposition
\ref{Natur} and Remark \ref{Trivial2}. The rest follows from Proposition
\ref{Trivial1}.
\end{proof}

\begin{corollary}
\label{Trivial4}Suppose $g_{*}=0$ in reduced homology. If $z\in H_{i}(X,A)$,
$i\neq n,$ then $\Lambda_{fg}(z)=0.$
\end{corollary}

\begin{proof}
It follows from Proposition \ref{Trivial1}.
\end{proof}

Thus, when $g$ is homologically trivial it suffices to consider only the
Lefschetz class with respect to an element $z\in H_{n}(X,A)$, see
\cite[Section 5]{Sav}. In fact the following condition is sufficient for
$\Lambda_{fg}$ to be nontrivial (cf. Proposition \ref{f-inverse}):

\begin{description}
\item [(A)]$f_{*}:H_{n}(X,A)\rightarrow H_{n}(M,\partial M)$ is a nonzero homomorphism.
\end{description}

\begin{proposition}
\label{cond(A)}\cite[Corollary 5.1]{Sav} If $f$ satisfies condition (A) and
$g_{*}=0$ in reduced homology then $\Lambda_{fg}\neq0$.
\end{proposition}

We call a map $f:(X,A)\rightarrow(M,\partial M)$ \textit{weakly
coincidence-producing} if every map $g:X\rightarrow M$ with $g_{*}=0$ has a
coincidence with $f$ (compare to coincidence producing maps \cite[Section
7]{BS}). Now we can restate Proposition \ref{cond(A)}.

\begin{corollary}
\label{CP}If $f$ satisfies condition (A) then $f$ is weakly coincidence-producing.
\end{corollary}

Let's consider some examples of applications of this corollary.

\begin{corollary}
\label{CPhmot}\cite[Corollary 5.6]{Sav} Suppose $M$ is a homotopy sphere,
$f:X\rightarrow M$ is a map, and

\begin{description}
\item [(A$^{\prime}$)]$f_{\#}:\pi_{n}(X)\rightarrow\pi_{n}(M)$ is onto$. $
\end{description}

\noindent Then condition (A) is satisfied, so $f$ is weakly coincidence-producing.
\end{corollary}

The proposition below follows from Lemma 5 of Gottlieb \cite{Gottlieb}.

\begin{corollary}
\label{CPsmooth}Let $M,X$ be smooth closed manifolds, $f:X\rightarrow M$ be
smooth, $\dim X=N,N>n.$ Suppose $F=f^{-1}(y)$ is a fiber with $y\in M$ a
regular point (then $F$ is a closed $(N-n)$-manifold) and

\begin{description}
\item [(A$^{\prime\prime}$)]$i^{\ast}:H^{N-n}(X)\rightarrow H^{N-n}(F)$ is nonzero.
\end{description}

\noindent Then condition (A) is satisfied, so $f$ is weakly coincidence-producing.
\end{corollary}

\begin{corollary}
\label{CPviet}Let $f:X\rightarrow M$ be an orientable fibration with fiber
$F.$ Suppose $F$ is arcwise connected, $A=f^{-1}(\partial M),$ $\partial
M\neq\emptyset,H_{i}(F)=0$ for $0<i<n-1.$ Then condition (A) is satisfied, so
$f$ is weakly coincidence-producing.
\end{corollary}

\begin{proof}
Suppose $\mathcal{C}$, a Serre class, contain only the zero group, see
\cite[Theorem 9.6.10, p. 506]{Spanier}. Then $\mathcal{C}$ is an ideal of
abelian groups. Observe that $H_{i}(M,\partial M)\in\mathcal{C}$ for $0\leq
i<1,$ and $H_{j}(F)\in\mathcal{C}$ for $0<j<n.$ Then $f_{\ast}:H_{q}%
(X,A)\rightarrow H_{q}(M,\partial M)$ is an $\mathcal{C}$-epimorphism for
$q\leq n$, so (A) is satisfied.
\end{proof}

\section{Examples.}

Brown \cite{Brown0} proved that a compact closed manifold $M$ is suitable (see
\cite{FN}) if and only if there is a multiplication on $M$ such that:

(1) $x\cdot e=x,\ \forall x\in M;$

(2) $\forall a,b\in M,\ \exists x\in M$ such that $a\cdot x=b$;

(3) $x\cdot y=x\cdot z\Rightarrow y=z,\ \forall x,y,z\in M.$

\noindent Then $M$ is an H-manifold and for any $x\in M$ there is a unique
$x^{-1}\in M$ such that $x\cdot x^{-1}=e$. The following is a slight
generalization of Theorem 3 of Wong \cite{Wong}.

\begin{theorem}
\label{Wong}If $M$ is a suitable manifold, $A=\varnothing,$ then
\[
\Lambda_{fg}(z)=<\overline{O_{M}},\psi_{\ast}(z)>,\ z\in H_{n}(X),
\]
where $\psi(x)=g(x)\cdot\lbrack f(x)]^{-1},\ x\in X,$ and $\overline{O_{M}} $
is the dual of the fundamental class $O_{M}.$
\end{theorem}

\begin{proof}
Consider the following commutative diagram:
\[%
\begin{array}
[c]{ccccccc}%
X\times X & ^{\underrightarrow{\ \ \ \ \ \ f\times g\ \ \ \ \ \ }} & M\times
M & ^{\underrightarrow{\ \ \ \ \ \ I\ \ \ \ \ \ }} & M^{\times} &
^{\ \underleftarrow{\ \ \ j\ \ \ }}\  & M,\\
\ \ \uparrow^{\delta} &  & \ \ \downarrow^{\sigma} &  & \ \ \downarrow
^{\sigma} & \swarrow_{k} & \\
X & ^{\underrightarrow{\ \ \ \ \ \ \psi\ \ \ \ \ \ }} & M & ^{\underrightarrow
{\ \ \ \ \ \ k\ \ \ \ \ \ }} & (M,M\backslash\{e\}) &  &
\end{array}
\]
where $\sigma(a,b)=b\cdot a^{-1},\ j(y)=(e,y),\ k$ is the inclusion. Since
$k_{\ast}:H_{n}(M)\rightarrow H_{n}(M,M\backslash\{e\})$ is an isomorphism,
for any $z\in H_{n}(X)$ we have the following
\[
I_{fg}(z)=I_{\ast}(f\times g)_{\ast}\delta_{\ast}(z)=j_{\ast}\psi_{\ast}(z).
\]
Therefore
\[
\Lambda_{fg}(z)=\varphi j_{\ast}\psi_{\ast}(z)=<\overline{O_{M}},\psi_{\ast
}(z)>.\
\]
\end{proof}

Corollaries \ref{CP} - \ref{CPviet} and Theorem \ref{Wong} can be used to
prove the existence of coincidences of maps between manifolds of different
dimensions. However we do not use the whole Lefschetz homomorphism, only its
part in dimension $n=\dim M$, i.e., $\Lambda_{fg}:H_{n}(X,A)\rightarrow
H_{0}(M)$ (in other words, we need only the Lefschetz number with respect to a
$z\in H_{n}(X,A)$ as in \cite{Sav}). This is also true for Example 2.4 in
\cite{GNO}: if $A:\mathbf{S}^{3}\times\mathbf{S}^{2}\rightarrow\mathbf{S}^{2}$
is the action given by regarding $\mathbf{S}^{2}$ as the homogeneous space
$\mathbf{S}^{3}/\mathbf{S}^{1}$ arising from the Hopf principal bundle
$\mathbf{S}^{1}\rightarrow\mathbf{S}^{3}\rightarrow\mathbf{S}^{2},$ then the
Knill trace is $0$ in all dimensions except $0$. This means that the only
nonzero part of $\Lambda_{pA}$ ($p$ is the projection) is the following:
$H_{2}(\mathbf{S}^{3}\times\mathbf{S}^{2})\rightarrow H_{0}(\mathbf{S}^{2}).$

To show that other values of the Lefschetz homomorphism may be important
consider Examples 2.3 and 7.2 in \cite{GNO}. In the notation of the present
paper, the first one states the following. If $\mu:G\times G\rightarrow G$ is
the multiplication of a compact Lie group then for $u\in H_{n+m}(G)$ we have
\[
\Lambda_{p\mu}(u\times O_{G})=\left\{
\begin{array}
[c]{ll}%
0\text{ \ \ \ \ \ } & \text{if }m<n\\
(-1)^{n}u & \text{if }m=n.
\end{array}
\right.
\]
The second example implies that if $Y$ is a manifold with $\dim Y=k>0$ and
$\chi(Y)\neq0$ then for a certain map $g:X=\mathbf{S}^{m}\times\mathbf{S}%
^{m}\times Y\rightarrow M=\mathbf{S}^{m}\times Y,$ the Lefschetz homomorphism
$\Lambda_{pg}$ is nontrivial in dimension $m+n$ $(n=\dim M=m+k) $:
$H_{m+n}(X)\rightarrow H_{m}(M)$. Thus we have an example of the Lefschetz
homomorphism with nontrivial values in dimensions other than $0$ or $n$.

Note that the statement in \cite{GNO} is that the maps $A,$ $\mu,$ and $g$
above and all maps homotopic to them have coincidences with the corresponding
projections $p.$ We have proved a little more than that: these maps have
coincidences with all maps homotopic to $p.$ To take this even further from
the setting of parametrized maps, we can consider the composition of the above
maps with another map. For example, suppose $H$ is a topological space,
$k:H\rightarrow G$ is a map. Let $f=p(k\times k),$ $g=\mu(k\times k):H\times
H\rightarrow G,$ so neither is the projection. Then by Theorem \ref{Natur} we
have for $u\in H_{n+m}(H\times H)$:
\[
\Lambda_{fg}(u\times k_{\ast}^{-1}(O_{G}))=\Lambda_{p\mu}(k_{\ast}(u)\times
O_{G})=\left\{
\begin{array}
[c]{ll}%
0\text{ \ \ \ \ \ } & \text{if }m<n\\
(-1)^{n}k_{\ast}(u) & \text{if }m=n.
\end{array}
\right.
\]
\medskip

In view of the Knill-like trace representation of the Lefschetz homomorphism
(Remark \ref{KnillTrace2}) we obtain the following.

\begin{proposition}
\label{2n}For $z\in H_{2n}(X,A)$, we have $\Lambda_{fg}(z)=g_{\ast}(f^{\ast
}(\overline{O_{M}})\frown z).$
\end{proposition}

It follows that if $X$ is a compact orientable $2n$-manifold, $f^{\ast n}%
\neq0$ and $\ker g_{\ast n}=0$ then $\Lambda_{fg}(O_{X})\neq0,$ where $O_{X}$
is the fundamental class of $X.$

\begin{center}
\textbf{Final Remarks.}
\end{center}

(1) The statement of Theorem \ref{EvalFormula} holds for open subsets of
$\mathbf{R}^{q}$, see \cite[Sections 6-10]{Sav}. As a result, the second part
of Theorem \ref{LefThm} can be proven for spaces more general than manifolds
(such as ANR's) by following Gorniewicz \cite[Sections V.2 and V.3]{Gorn}. It
would also be interesting to try to obtain Theorem \ref{LefThm} for a
non-orientable manifold $M$ by following Gon\c{c}alves and Jezierski \cite{GJ}.

(2) We know that the Hopf map $h:\mathbf{S}^{3}\rightarrow\mathbf{S}^{2}$ is
onto, in other words, it has a coincidence with any constant map $c$. On the
other hand, as $\Lambda_{fg}=0$ for any pair of maps $f,g:\mathbf{S}%
^{N}\rightarrow\mathbf{S}^{n}$ and $N\neq n$, the Lefschetz-type Coincidence
Theorem \ref{LefThm} fails to predict the existence of coincidences of
$(h,c)$. In fact, $h$ has a coincidence with any map homotopic to $c$
\cite{Brooks}, therefore the converse of the Lefschetz coincidence theorem for
spaces of different dimensions fails.

\end{document}